\documentclass[runningheads,a4paper]{llncs}
\usepackage{amssymb}
\usepackage{amsmath}
\usepackage{graphicx}
\usepackage{url}

\usepackage[T1]{fontenc} 
\usepackage[utf8]{inputenc}
\usepackage{color}

\begin{document}

\mainmatter 

\title{Nonlinear Feedback Shift Registers and Zech Logarithms}

\author{Janusz Szmidt}

\institute{ Military Communication Institute \\
ul. Warszawska 22A, 05-130 Zegrze Południowe, Poland \\
j.szmidt@wil.waw.pl} 

\maketitle

\begin{abstract}
We construct feedback functions of Nonlinear Feedback Shift Registers from those of Linear Feedback Shift Registers using the cross-join pairs method and the Zech logarithms in finite fields. We present a hypothetical algorithm to generate all de Bruijn sequences of a given order. 
\end{abstract}

\section{Introduction}
Recently several papers has apeared on construction of de Bruijn sequences from some classes of Linear Feedback Shift Registers (\textit{LFSRs}), see \cite{CEz,CRV,DP,LJL,LZHLH,LZLH,LZLHL,MS,SD}. The authors have used the method of joining cycles generated by some reducible \textit{LFSRs.} In \cite{MS} we have proved that having one de Bruijn sequence of a given order one can construct all de Bruijn sequences of that order by repeated application of the cross-join pairs method. We present a method to find the 
cross-join pairs of states by calculating the Zech logarithms in Galois fields. This enables us to construct feedback functions of Nonlinear Feedback Shift Registers (\textit{NFSRs}) of maximum period up to the order of $ n = 430 $ and for any primitive trinomial. The Freyers formula is a  new analytic tool in the theory and its consequences are presented.

\section{Basic definitions and results}
Let $ \mathbb{F}_{2} = \{0, 1 \} $ denote the binary field and $ \mathbb{F}_{2}^{n} $ the vector space of all binary $ n $-tuples.
A binary Feedback Shift Register (\emph{FSR}) of order $ n $ is a mapping
$\mathfrak{F} : \mathbb{F}_{2}^{n} \longrightarrow \mathbb{F}_{2}^{n} $
of the form
\begin{equation}
\mathfrak{F} : (x_{0}, x_{1}, \dots , x_{n-1}) \longmapsto (x_{1}, x_{2}, \ldots , x_{n-1}, F(x_{0}, x_{1}, \dots , x_{n-1})),
\end{equation}
where the \emph{feedback function} $ F $ is a Boolean function of $ n $ variables.
The \emph{FSR} is called \emph{non-singular} if the mapping $ \mathfrak{F} $  
is a bijection of $ \mathbb{F}_{2}^{n} $. 
It was proved in \cite{Gol}  that the non-singular feedback function of (1) has the form
\[F(x_{0}, x_{1}, \dots, x_{n-1}) = x_{0} + f(x_{1},\dots,x_{n-1}).\]

\begin{definition} 
A \emph{de Bruijn sequence} of order $ n $ is a sequence of length $ 2^{n} $ of elements of $\mathbb{F}_{2}$ in which all different $ n $-tuples appear exactly once. 
\end{definition}

\begin{definition} 
A \emph{modified} de Bruijn sequence of order $ n $ is a sequence of length$  $ $ 2^{n}-1 $ obtained from a de Bruijn sequence of order $ n $ by removing one zero from the tuple of $ n $ consecutive zeros. 
\end{definition}
It appears \cite{Gol} that de Bruijn sequences are generated by non-singular \textit{NFSRs} of a special type.
\begin{theorem}
Let $(s_t)$ be a de Bruijn sequence of order $ n $. Then there exists a  Boolean function \\
$f(x_1,\dots,x_{n-1})$ such that 
\[    s_{t+n} = s_t + f(s_{t+1},\dots,s_{t+n-1}),\ \ t=0,1, \dots,2^{n}-n-1.\]
\end{theorem}

It was proved by Flye Sainte-Marie \cite{F} in 1894 and independently by \mbox{de Bruijn} \cite{B} in 1946 that the number of cyclically non-equivalent sequences satisfying Definition 1 is equal to
\[B_{n} = 2^{2^{n-1}-n}.\] 

Let $(a_{t}) = (a_{0}, a_{1}, \dots, a_{2^n-1})$ be a de Bruijn sequence.
Let us denote $S_{i}=(a_{i}, a_{i+1}, \dots, a_{i+n-1})$. Then the de Bruijn sequence can be represented as the sequence of states
$(S_{t}) = (S_{0}, S_{1}, \dots, S_{2^{n}-2}, S_{2^{n}-1})$. 

Let $ \alpha = (a_{0},a_{1}, \dots , a_{n-1})$ be a state (i.e., $ n $ consecutive elements) of a de Bruijn sequence generated by the feedback function $ F. $ The conjugate of the state $ \alpha $ is $\widehat{\alpha} = (\overline{a}_0, a_1, \dots, a_{n-1}),$
where $ \overline{a} = a + 1 $ is the negation of the bit $ a. $

\begin{definition}
Two pairs of conjugate states of \emph{FSR} determined by the feedback function $ F : $ 
\[ \alpha = (a_0, a_1, \dots, a_{n-1}), \ \ \widehat{\alpha} = (\overline{a}_0, a_1, \dots, a_{n-1}), \] 
\[ \beta = (b_0, b_1, \dots, b_{n-1}), \ \ \widehat{\beta} = (\overline{b}_0, b_1, \dots, b_{n-1}), \]
constitute \emph{cross-join pairs} for the sequence $ (a_{t}) $ generated by the \emph{FSR} if the states occur in the order
$ \alpha, \beta,  \widehat{\alpha}, \widehat{\beta} $ during the generation of the sequence.
\end{definition}

For a given de Bruijn sequence (or a modified de Bruijn sequence) and its cross-join pairs of states one can construct a new de Bruijn (modified de Bruijn) sequence by applying the cross-join pairs operation \cite{Gol}. Mykkeltveit and Szmidt \cite{MS} have proved the following 

\begin{theorem}
Let $(u_t),(v_t)$ be two different de Bruijn sequences of order $n$. Then  $(v_t)$  can be obtained from  $(u_t)$  by repeated application of the cross-join pairs operation.
\end{theorem}

\section{Our contribution}

Let $p(x) = x^{n} + c_{n-1}x^{n-1} + \dots + c_{1}x + 1 $ \
be a primitive polynomial of degree $ n $ with binary coefficients.
Then the linear recursion
\[g(x_{0}, x_{1}, \ldots , x_{n-1}) = x_{0} + c_{1}x_{1} + \dots + c_{n-1}x_{n-1} \]
generates the maximum period sequence, it is called an $m$-sequence.
Let $ a $ be a root of the polynomial $ p(x) $, i.e. $ p(a) = 0 $ in the Galois field $GF(2^n).$
The sequence of elements of the field $ GF(2^n): $
\[ \{1,a.a^2, \ldots , a^{2^{n}-2} \}  \]
has period $ 2^n-1 $ and from it the binary $ m $-sequence can be constructed.

For an integer $ j \in \{1, \ldots , 2^n-2 \} $ we define the integer $ Z(j) $ such that
$ 1 + a^j = a^{Z(j)} $. Then we have a one-to-one function
\[  Z : \{1, \ldots , 2^n-2 \}  \longrightarrow     \{1, \ldots , 2^n-2 \} \]
which is called the Zech logarithm or Jacobi logarithm.
The Zech logarithms are tabularized. There are effective algorithms to calculate them. Magma can be used to calculate the Zech logarithms for the order $ n $ up to 430, i.e. in the field $GF(2^{430})$.

We show that working knowledge of Zech logarithms enables us to find cross-join pairs of states for a primitive \textit{LFSR} and then to construct feedback functions of \textit{NFSRs} generating modified de Bruijn sequences. We give two examples of this construction.

Let us take the primitive polynomial $ p(x) = x^{31} + x^3 +1$ and its root $p(a) = 0.$
We have the following cross-join pairs: $ c = (3,6,31,62),$ which is an abbreviation for two pairs of states $(a^{3}, a^{6}, 1+a^{3}=a^{31}, 1+a^{6}=a^{62}) $ since 
$ Z(3) = 31, \ Z(6) = 62, $ where we have used the formula $ Z(2n) = 2Z(n) $ for the Zech logarithm.
The states of \textit{LFSR} at 'jumps' where the linearity is broken are
\[ A = (0, 0, 0, 0, 0, 0, 0, 0, 0, 0, 0, 0, 0, 0, 0, 0, 0, 0, 0, 0, 
0, 0, 0, 0, 0, 0, 0, 0, 1, 0, 0), \ A_{28} = 1, \]
\[ B =  (0, 0, 0, 0, 0, 0, 0, 0, 0, 0, 0, 0, 0, 0, 0, 0, 0, 0, 0, 0, 
0, 0, 0, 0, 0, 1, 0, 0, 0, 0, 0), \ B_{25} = 1. \]

The corresponding feedback function of the constructed \textit{NFSR} is
\[ F = x_0 + x_3 +  \prod^{30}_{i=1} (x_i + A_i + 1) + \prod^{30}_{i=1} (x_i + B_i + 1).\]
It is a Boolean function of algebraic degree 29.

In the second example we construct a sequence of non-overlapping cross-join pairs. We take
the primitive polynomial $ p(x) = x^{127} + x + 1. $
Since $ Z(1) = 127, \ Z(2) = 254, $ we have the sequence of cross-join pairs: 
\[ c_i = (2^{8i}, 2^{1+8i}, 127 \cdot 2^{8i}, 127 \cdot 2^{1+8i}), \ \ i = 0,1, \ldots , 15,\]
which are mutually disjoint.
From this family we can construct $ 2^{16}-1 $ \textit{NFSRs} of order $ n = 127 $ which generate sequences of period $ 2^{127}-1. $
As an example we have the cross-join pairs $ c_3 = (2^{24},2^{25},127 \cdot 2^{24}, 127 \cdot 2^{127}).$ The corresponding Boolean feedback function has algebraic degree 125.

\section{Fryers' Formula}
Let $ S(n) $ be the set of functions $ \mathbb{F}_{2}^{n-1} \rightarrow \mathbb{F}_{2}  $
that generate a de Bruijn sequence of order $n$. For $f : \mathbb{F}_{2}^{n-1} \rightarrow \mathbb{F}_{2} $, define $ S(f;k) $ to be the set of $ g \in S(n) $ such that the weight of $ f + g $
is $ k $. In other words, the number of disagreements between the functions $ f $ and $ g $
is $ k $. Moreover, let $ N(f;k) = |S(f;k)| $ and 
\begin{equation}
 G(f;y) = \sum_{k} N(f;k)y^k .
\end{equation} 
Let a linear function $ l : \mathbb{F}_{2}^{n-1} \rightarrow \mathbb{F}_{2} $ generate a
$(2^n - 1) $ - periodic sequence without a run of $n$ zeros.

In \cite{CRV} a formula due to Michael Fryers is given: 
\begin{equation}
G(l;y) = \frac{1}{2^n} \left((1 + y)^{2^{n-1}} - (1 - y)^{2^{n-1}} \right)
\end{equation}
which can be deduced from Theorem 1.1 of \cite{CRV}.
The authors of \cite{CRV} view it as a generalization of the result of Helleseth and Kl\o ve \cite{HK} who computed $ N(l;3) $
for any such $ l. $
We derive a formula for $ N(l;k) $ using the Maclaurin expansion of $ G(l;y): $
\begin{equation}
G(l;y) = \sum_{k = 1}^{\infty} a_k y^k , \ \ \ \ a_k = \frac{1}{k!}G^{(k)}(l;y)\mid_{(y=0)} .
\end{equation}
Hence $ N(l;k) = a_k $ can be obtained by calculating the derivatives of $ G(l;y) $ at $ y=0. $
We have: \\
\begin{enumerate}
\item $ N(l;1) = 1, $  \ \ which means that from an $ m $-sequence we get one de Bruijn sequence by adding the all zero state which makes one change in the truth table of the function $ l $.
\item $ N(l;2) = 0; $ \ \ now we change the truth table in two places, one adds the zero cycle and the other disjoins the whole cycle into two cycles. In general, $ N(l;k) = 0 $ for all even $k$, since we obtain some number of disjoints cycles.
\item \[ N(l;3) = \frac{(2^{n-1} - 1)(2^{n-1} - 2)}{3!}\] 
which is the Helleseth and Kl\o ve formula \cite{HK} which counts the number of cross-join pairs in an $ m $-sequence and equivalently the number of double changes in the truth table (plus one which leads to the zero cycle) which disjoin the cycle and then join it to one cycle. Hence the Fryers formula (3) leads to an analytic proof of the Helleseth and Kl\o ve formula.
\item \[ N(l;5) = \frac{(2^{n-1} - 1)(2^{n-1} - 2)(2^{n-1} - 3)(2^{n-1} - 4)}{5!} \]
which is a formula for an application of the cross-join pairs method twice and the number of de Bruijn sequences obtained this way.
\item \[ N(l,k) = \frac{(2^{n-1} - 1)(2^{n-1} - 2)(2^{n-1} - 3)(2^{n-1} - 4)\dots (2^{n-1} - k + 2)(2^{n-1} - k + 1)}{k!} \]
for $ k = 3,5. \dots , 2^{n-1} - 1.$ It is a formula for the number of  de Bruijn sequences obtained  after application of the cross-join pairs method $ (k-1)/2 $ times.
\item $ N(l;2^{n-1} - 1) = 1, $ which follows from a calculation; it means that at the last step
we get the last new de Bruijn sequence.
\item $ N(l,k) = 0, $ for $ k >  2^{n-1} - 1, $ which means that no more de Bruijn sequences are 
produced.
\item We have
\begin{equation}
G(l,1) = \sum_{k=1}^{2^{n-1}-1} N(l,k) = 2^{2^{n-1} - n},
\end{equation}
which is the number of all cyclically non-equivalent de Bruijn sequences of order $ n. $
\end{enumerate}

\subsection{The special cases}
Let us consider some special cases of the equation (5) for small values of $ n. $ \\
\begin{enumerate}
\item $ n = 4 : $
\[  \sum_{k=1}^{7} N(l,k) = 1 + 7 + 7 + 1 = 16. \]
\item $ n = 5 : $
\[ \sum_{k=1}^{15} N(l,k) = 1 + 35 + 273 + 715 + 715 + 273 + 35 + 1 = 2048.\]
Here 35 is the number of cross-join pairs for an $ m $-sequence of period $ 2^5 - 1. $ 273 is the number of new de Bruijn sequences obtained after second application of the cross-join pairs method. Application of the cross-join pairs method seven times produces all de Bruijn sequences of order 5.
\item $ n = 6 : $
\[ \sum_{k=1}^{31} N(l,k) = 1 + 155 + 6293 + 105183 + 876525 + 4032015 + 10855425 + 17678835 + \] 
\[ \ \ \ \ \ \ \ \ \ \ \ \ \ \ \ \ 17678835 + 10855425 + 4032015 +  876525 + 105183 + 6293 + 155 + 1\]
\ \ \ \ \ \ \ \ \ \ \ \ \ \ \ \ \ \ \ \ \ \ \ \ \ \ $ = 67108864 = 2^{26}.$ \\
Observe a symmetry in the appearing numbers.
\end{enumerate}
We have confirmed experimentally the appearance of the above numbers for $ n = 4, $ and of the number 273  for $ n = 5 $. The consequences of the Freyers formula are compatible with the theorem proved in \cite{MS} which states that starting from one de Bruijn sequence of a given order one can obtain all de Bruijn sequences of that order by applying the cross-join pairs method several times. 

\subsection{ A general Algorithm}

\begin{enumerate}
\item Choose an $ m $-sequence of a given order $ n $ and the corresponding linear recursion $ l. $ Find all cross-join pairs for the corresponding modified de Bruijn sequence (of period $2^n - 1$).
Construct all Boolean feedback fuctions for those cross-join pairs. There are $ N(l,3) $ of them.
\item Find all cross-join pairs for the just constructed modified de Bruijn sequences and obtain feedback functions from those cross-join pairs. It happens that some \emph{NLFSRs} produced by different cross-join pairs are identical. Hence we have to sieve the set of all \emph{NLFSRs} obtained in second application of the cross-join pairs method to have only different feedback functions. Their number should be $ N(l;5). $
\item Realize the next steps of application of cross-join pairs method up to the last one where only one new \emph{NLFSRs} should appear. This way all feedback functions of modified de Bruijn sequences of order $ n $ have to be generated.
\end{enumerate}

This algorithm is a hypothetical one since practically we can generate all feedback functions of de Bruijn sequences for $ n = 5 $ and possibly for $ n = 6. $ For $ n = 7 $ it is impractical since there are $ 2^{57} $ de Bruijn sequences of order 7. One can try to follow this algorithm for greater values of $ n $ and possibly for greater values of $ k $ by applying the cross-join pairs method only for some subsets of cross-join pairs. A theoretical and practical problem appears how to construct cross-join pairs for level $ k. $
 
\section{Conclusions}
We construct feedback functions of \textit{NFSRs} from those of \textit{LFSRs} by applying the cross-join pairs method. Our contribution is using the Zech logarithms to find the cross-join pairs of states of \textit{LFSRs} which are described by primitive polynomials. This enables us to construct feedback functions in the range of orders of \textit{NFSRs} where the Zech logarithms can be calculated and for all known primitive trinomials. We present some consequences of the Fryers formula.

\end{document}